\documentstyle[11pt]{article}
\begin{document}

\begin{titlepage}

\hfill{UTAS-PHYS-96-08}\\
\mbox{}\hfill{}
\vskip 1.6in
\begin{center}
{\Large {\bf On boson algebras as Hopf algebras }}\\[5pt]
{\Large {\bf }}
\end{center}

\normalsize
\vskip .4in

\begin{center}
I Tsohantjis$^{\dag}$\hspace{3pt}
A Paolucci$^{*}$  \hspace{3pt}
and \hspace{3pt} P D Jarvis$^{\dag}$
\par \vskip .1in \noindent
$^{\dag}${\it Department of Physics, University of Tasmania}\\
{\it GPO Box 252C Hobart, Australia 7001}\\
$^{*}${\it School of Mathematics, University of Leeds,}\\
{\it Leeds, LS2 9JT UK}
\end{center}
\par \vskip .3in

\begin{center}
{\Large {\bf Abstract}}\\
Certain types of generalized undeformed and deformed boson algebras  
which admit a Hopf algebra structure are introduced, together with their 
Fock-type representations and their corresponding $R$-matrices. It is also shown that a class of 
generalized Heisenberg 
algebras including those underlying
physical models such as that of Calogero-Sutherland, is isomorphic with one 
of the types of boson algebra proposed, and can be formulated as a Hopf algebra.

\end{center}

\vspace{1cm}

\begin{center}
{\bf Physics and Astronomy Classifications}. 02.10.
\end{center}

\end{titlepage}

\section{Introduction.}

Deformations of the boson algebra have been recently the subject of
extensive research partly because of their significance in quantum groups
(see for example \cite{drin}-\cite{book}) and supergroups \cite{gould}.
Chronologically first comes the Arik-Coon $q$-deformation of the Heisenberg
algebra \cite{aric}:

\begin{equation}
aa^{\dagger }-qa^{\dagger }a=I,  \label{arik}
\end{equation}
followed by the Macfarlane-Biedenharn \cite{macf,bid}, and Sun and Fu \cite
{sun} $q$-deformed bosons

\begin{equation}
aa^{\dagger }-q^{\pm }a^{\dagger }a=q^{\mp N},  \label{sunfu}
\end{equation}
the Chakrabarti-Jagannathan two-parameter model \cite{chak}

\begin{equation}
aa^{\dagger }-pa^{\dagger }a=q^{-N},
\end{equation}
and the Calogero-Vasiliev model \cite{vas1}

\begin{equation}
\lbrack a,a^{\dagger }]=I+2\nu K,\quad K^2=I,
\end{equation}
which coupled with (\ref{sunfu}) as

\begin{equation}
aa^{\dagger }-qa^{\dagger }a=q^{-N}(I+2\nu K),
\end{equation}
was studied in \cite{mac2}, while its $q$-deformation by Macfarlane \cite
{macfar} is

\begin{equation}
aa^{\dagger }-q^{\pm (I+2\nu K)}a^{\dagger }a=[I+2\nu K]_qq^{\mp (N+\nu -\nu
K)},\quad K=(-1)^N,
\end{equation}
where as usual $[x]_q=(q^x-q^{-x})/(q-q^{-1})$. In addition with the
Katriel-Quesne minimally deformed oscillators \cite{que} which provides an
attempt to unify existing deformed oscillators, generalizations and
applications of the above models to mathematics and physics have been of
increasing interest, and their consistency, interrelation and
representations have been well analyzed \cite{kur}, \cite{jan} \cite{ng}-%
\cite{solo}.

Generalizations of the usual Heisenberg algebra that have appeared in \cite
{vas1} have been implemented \cite{pl1}, \cite{pl2}, in describing
relativistic fields with arbitrary fractional spin (anyon), ``bosonizing''
supersymmetric quantum mechanics and pointing out the relation of it with
integrable quantum mechanical models \cite{brink} such as the
Calogero-Sutherland model \cite{cal1}, \cite{cal2}.

On the other hand, the recent investigation \cite{pal1}-\cite{pal5} of
simplest $q$-deformations of the Heisenberg algebra has also been shown to
play a key part in obtaining and classifying representations of deformed
bosons algebras \cite{chfelp}.

The relation of a possible Hopf algebra structure consistent both with an
appropriate definition of a boson algebra and its deformation have also been
addressed \cite{oh}, \cite{palev1}-\cite{yan3}. In particular in \cite{yan1}%
, \cite{yan2}, \cite{yan3} a certain definition of deformed boson algebra
was investigated having a Hopf algebra structure, while in \cite{dt} the $R$%
-matrix obtained was corrected and generalized using the quantum double
construction. The results of \cite{yan1} were partly generalized in \cite{oh}%
.

The aim of this paper is to investigate certain generalizations of
undeformed and deformed boson algebras possessing a Hopf algebra structure,
which in \cite{pa} will be used to establish an algebraic relation with
already known boson algebras (undeformed and deformed). In section two of
this paper we introduce general notions on quasitriangular Hopf algebras. In
section three we introduce and analyze the properties of undeformed
generalized boson algebras $B_{\alpha ,\beta }$ and $\overline{B}_{\sigma
,\tau }$ which admit a Hopf algebra structure, while in section four we give
a $q$-deformation of the previous algebras, prove that they also admit a
Hopf structure and present an $R$-matrix for each one of them. We further
analyze in section five a more general form of the `deformed' Heisenberg
algebra $H_\nu $ of \cite{vas1}, showing that under certain conditions it
admits a Hopf algebra structure and demonstrate its connection with the
undeformed boson algebra $B_{\alpha ,\beta }$ defined in the third section.
Finally we end with certain comments on possible physical and mathematical
applications and consequences of our approach.

\section{Generalities on quasitriangular Hopf algebras.}

Consider a unital associative algebra, over a field $F$ , with
multiplication $m:A\otimes A\rightarrow A$ (i.e. $m(a\otimes b)=ab$, $%
\forall a$ and $\forall b$ $\in A$) and unit $u:F\rightarrow A$ (i.e. $%
u(1)=I $, the identity on $A$) endowed with a Hopf algebra structure, that
is, having a coproduct $\Delta :A\rightarrow A\otimes A$, a counit $%
\varepsilon $ :$A\rightarrow F$ (which is a homomorphism) and an antipode $S$%
: $A\rightarrow A$ (which is an antihomomorhism i.e. $S(ab)=S(b)S(a)$, and
we shall assume that it has an inverse $S^{-1}$) subject to the following
consistency condition:

\[
(id\otimes \Delta )\Delta (a)=(\Delta \otimes id)\Delta (a) 
\]

\[
(id\otimes \varepsilon )\Delta (a)=(\varepsilon \otimes id)\Delta (a)=a 
\]

\begin{equation}
m(id\otimes S)\Delta (a)=m(S\otimes id)\Delta (a)=\varepsilon
(a)I \quad \forall a\in A.  \label{hop}
\end{equation}
Following Sweedler\cite{swed} we write

\[
\Delta (a)=\sum_{(a)}a^{(1)}\otimes a^{(2)} 
\]

\begin{equation}
\Delta _n(a)=(\Delta \otimes I^{\otimes (n-1)})\Delta
_{n-1}(a)=\sum_{(a)}a^{(1)}\otimes a^{(2)}\ldots \otimes a^{(n+1)}
\label{cop}
\end{equation}
Let $T$ be the twist map on $A\otimes A$ defined by $T(a\otimes b)=b\otimes
a $. Then there also exists an opposite Hopf algebra structure on $A$ with
coproduct $T\Delta \doteq \Delta ^T$ , antipode S$^{-1}$ and counit as
before. According to Drinfeld \cite{drin} a Hopf algebra A is called
quasitriangular if there exists an invertible element $R$ such that

\begin{equation}
R=\sum_ia_i\otimes b_i\in A\otimes A  \label{r1}
\end{equation}

\begin{equation}
\Delta ^T(a)R=R\Delta (a)\quad\forall a\in A  \label{r2}
\end{equation}

\begin{equation}
R_{12}R_{13}R_{23}=R_{23}R_{13}R_{12}  \label{r3}
\end{equation}
where $a_i$ and $b_i$ are appropriately chosen dual basis of $A$, $%
R_{13}R_{23}=(\Delta \otimes I)R$, $R_{13}R_{12}=(I\otimes \Delta )R$ and $%
R^{-1}=(S\otimes I)R$. Finally for any Hopf algebra we can define the
adjoint operations $ad$ and $ad^{\prime }$ given by

\begin{equation}
ad_a(b)=\sum_{(a)}a^{(1)}bS(a^{(2)})\quad  \label{ad}
\end{equation}

\begin{equation}
ad_a^{\prime }(b)=\sum_{(a)}a^{(2)}bS^{-1}(a^{(1)})\quad \forall a\in A.
\label{ad'}
\end{equation}

\section{The undeformed generalized boson algebras $B_{\alpha ,\beta }$ and $%
\overline{B}_{\sigma ,\tau }$ .}

We start by recalling that the boson algebra $B$ is generated by $a$, $%
a^{\dagger }$, and $N$ subject to the following relations:

\begin{equation}
\lbrack a,a^{\dagger }]=I,\quad [N,a]=-a,\quad [N,a^{\dagger }]=a^{\dagger }
\label{normal}
\end{equation}
and a Fock space representation is provided by

\begin{equation}
|n>=\frac 1{\sqrt{n!}}(a^{\dagger })^n|0>,\quad N|n>=n|n>,\quad a|n>=\sqrt{n}%
|n-1>,\quad a^{\dagger }|n>=\sqrt{(n+1)}|n+1>.  \label{focnormal}
\end{equation}

The often used identification $a^{\dagger }a=N$ and $aa^{\dagger }=N+1$
holds in the quotient $B/<C>$ (and on the above Fock space) where $<C>$ is
the two-sided ideal generated by $C=a^{\dagger }a-N$. As was demonstrated in 
\cite{palev1}, a Hopf algebra structure on this algebra fails to exist.

\subsection{The algebra $B_{\alpha ,\beta }$.}

We shall now consider the family of algebras $B_{\alpha ,\beta }$ generated
by $a,$ $a^{\dagger }$ and $N$ subject to the following relations: 
\begin{eqnarray}
\{a,a^{\dagger }\} &=&\alpha N+\beta I,  \nonumber  \label{bos} \\
\lbrack N,a] &=&-a,  \label{bos1} \\
\lbrack N,a^{\dagger }] &=&a^{\dagger }  \nonumber
\end{eqnarray}
where $\alpha ,$ $\beta \in {\bf R}$ . Here and in the rest of the paper $%
\{x,y\}=xy+yx$ . If we take the quotient of $B_{2,1}$ with respect to an
ideal generated by $a^{\dagger }a-N$ we recover $B/<C>$ above. This algebra
can be enlarged to become a Hopf algebra by adding an invertible element $%
(-1)^N$ which will be treated as a supplementary generator satisfying the
following relations:

\begin{eqnarray}
\{(-1)^N,a\}=0=\{(-1)^N,a^{\dagger }\},\quad [(-1)^N,N]=0,\quad  \label{n}
\end{eqnarray}
Similar considerations were used in \cite{maj} and in that paper's context
our enlarged algebra $B_{\alpha ,\beta }^{}$ can be thought of as a spectrum
generating algebra for the ordinary harmonic oscillator, while the element $%
g $ of \cite{maj} will be $g=(-1)^{\widetilde{N}}$ provided that we impose
the condition $g^2=(-1)^{2\widetilde{N}}=I$ where $\widetilde{N}=N+\frac
\beta \alpha $ . At this point we do not necessarily have to impose this
condition (which implies that $(-1)^{2N}=(-1)^{-\frac{2\beta }\alpha }).$ We
shall denote by $B_{\alpha ,\beta }^{+}$ and $U(B_{\alpha ,\beta }^{+})$,
this enlarged algebra and its universal enveloping algebras respectively.
Then the coproduct counit and antipode satisfy (\ref{hop}) are given by:

\[
\Delta (N)=N\otimes I+I\otimes N+\frac \beta \alpha I\otimes I 
\]

\[
\Delta (a)=a\otimes I+(-1)^{N+\frac \beta \alpha }\otimes a 
\]

\[
\Delta (a^{\dagger })=a^{\dagger }\otimes I+(-1)^{-N-\frac \beta \alpha
}\otimes a^{\dagger } 
\]

\[
\varepsilon (N)=-\frac \beta \alpha ,\quad \varepsilon (a)=\varepsilon
(a^{\dagger })=0,\quad \varepsilon (I)=1 
\]

\begin{equation}
S(N)=-N-\frac{2\beta }\alpha ,\quad S(a)=-(-1)^{-N-\frac \beta \alpha
}a,\quad S(a^{\dagger })=a^{\dagger }(-1)^{N+\frac \beta \alpha },
\label{hop1}
\end{equation}

\begin{equation}
\Delta ((-1)^{\pm \widetilde{N}})=(-1)^{\pm \widetilde{N}}\otimes (-1)^{\pm 
\widetilde{N}},\quad \varepsilon ((-1)^{\pm \widetilde{N}})=I,\quad
S((-1)^{\pm \widetilde{N}})=(-1)^{\mp \widetilde{N}}.  \label{n'}
\end{equation}
provided that $\alpha \neq 0.$ Moreover an opposite Hopf algebra structure
also exists for $B_{\alpha ,\beta }$ with coproduct $\Delta ^T$ and antipode
the inverse $S^{-1}$of $S$ given by :

\begin{eqnarray}
S^{-1}(N) &=&S(N),\quad S^{-1}(a)=(-1)^{-N-\frac \beta \alpha }a,\quad 
\nonumber \\
S^{-1}(a^{\dagger }) &=&-a^{\dagger }(-1)^{^{N+\frac \beta \alpha }},\quad
S^{-1}((-1)^{\pm \widetilde{N}})=S((-1)^{\pm \widetilde{N}}).
\end{eqnarray}
A Fock-type representation, with $a|0>=0,$ $|n>^{\prime }=(a^{\dagger
})^n|0>,$ $N|n>^{\prime }=n|n>^{\prime },n\in Z_{+}$ and $<0|0>=1$, exists
such that,

\begin{equation}
a|n>^{\prime }=[n]|n-1>^{\prime },\quad {\mbox{\rm where}} 
\quad [n]=\frac{\alpha n}
2+\frac{2\beta -\alpha }4(1+(-1)^{n+1})  \label{foc1'}
\end{equation}
If $\alpha >0,$ $\beta >0$ a unitary representation of $B_{\alpha ,\beta }$
is provided by:

\[
|n>=\frac 1{([n]!)^{\frac 12}}(a^{\dagger })^n|0> 
\]

\[
a|n>=[n]^{\frac 12}|n-1>,\quad a^{\dagger }|n>=[n+1]^{\frac 12}|n+1>, 
\]

\begin{equation}
{\mbox{\rm where}} \quad [n]!=\prod_{l=1}^n[l]\quad {\mbox{\rm and}}\quad <n|n^{\prime
}>=\delta _{nn^{\prime }}.  \label{foc1}
\end{equation}
With the definition $(-1)^{\pm N}|n>=(-1)^{\pm n}|n>$ this Fock space
provides also a representation of $B_{\alpha ,\beta }^{+}$.

There exists an element $L$ in the enveloping algebra $U(B_{\alpha ,\beta })$
of $B_{\alpha ,\beta }$ given by

\begin{equation}
L=\lambda _1a^{\dagger }a+\lambda _2N+\lambda _3I,  \label{L}
\end{equation}
and such that

\begin{equation}
\{L,a\}=\{L,a^{\dagger }\}=0,  \label{La}
\end{equation}
provided that the following constraints on $\lambda _i\in {\bf R}$ $%
(i=1,2,3) $ are satisfied

\begin{equation}
2\lambda _3+\beta \lambda _1+\lambda _2=0,\quad 2\lambda _2+\alpha \lambda
_1=0,  \label{con}
\end{equation}
(these will give for example $\lambda _2/\lambda _1=-\alpha /2,$ and $%
\lambda _3/\lambda _1=(\alpha -2\beta )/4=-$ $\alpha /2$ ($\lambda
_3/\lambda _2),$ with $\alpha \neq 0)$. This choice of $L$ subject to (\ref
{con}) is obviously not unique as it can be easily checked that any odd
power of $L$ will satisfy (\ref{La}). However (\ref{L}) is the {\it unique}
element of linear combination of lowest order monomials of generators of $%
B_{\alpha ,\beta }$ that will satisfy (\ref{La}). This can be inferred by
writing $L^{\prime }=C_{lmn}K^l(a^{\dagger })^ma^n$ , $l$, $m$, $n$ $\in 
{\bf Z}_{+}$ , $C_{lmn}\in {\bf R}$ and demanding that (\ref{La}) are
satisfied together with $[L^{\prime },N]=0$. For a given $B_{\alpha ,\beta
}, $ i.e. for given values of $\alpha $ and $\beta $ relations (\ref{con})
give us the conditions on $\lambda _i$ under which $L$ becomes zero. On what
follows we shall assume unless otherwise stated that $L$ is non-zero (e.g.
when $\alpha =0$ and $\beta =0$ then $L\neq 0$ if and only if $\lambda
_1\neq 0$ or when $\alpha =0$ and $\beta \neq 0$ then $L\neq 0$ if and only
if $\beta \lambda _1=-2\lambda _3\neq 0$). Then using (\ref{con}), $L\in
U(B_{\alpha ,\beta })$ can be put in the form

\begin{equation}
L=\lambda _1\left( a^{\dagger }a-\frac \alpha 2N+(\frac \alpha 4-\frac \beta
2)I\right) ,\quad \lambda _1\neq 0.  \label{L'}
\end{equation}
If we consider $B_{\alpha ,\beta }^{+}$ then the additional term $\lambda
_4(-1)^{\widetilde{N}}$ ($\lambda _4$ $\in {\bf R},\lambda _4\neq 0$) can be
considered and an element $L^{+}\in U(B_{\alpha ,\beta }^{+})$ can be taken
as

\begin{equation}
L^{+}=L+\lambda _4(-1)^{\widetilde{N}}  \label{L''}
\end{equation}
satisfying (\ref{La}), while for given values of $\alpha $ and $\beta ,$ $%
L^{+}$ is also not unique and odd powers of it will give (\ref{La}). However
it should be noted that in the case of $B_{\alpha ,\beta }^{+}$ the element $%
\lambda _4(-1)^{\widetilde{N}}$ is the unique non-zero lowest order monomial
satisfying (\ref{La}). This again can be inferred by writing $L^{\prime
+}=C_{plmn}(-1)^{pN}K^l(a^{\dagger })^ma^n$, $p$, $l$, $m$, $n$ $\in {\bf Z}%
_{+}$, $C_{plmn}\in {\bf R}$ and demanding that (\ref{La}) are satisfied
together with $[L^{\prime +},N]=0$. Relation (\ref{L''}) is then the next
most general one to be considered. Utilizing (\ref{hop1}) and (\ref{con}),
it can be found that $\Delta (L^{+}),$ $S(L^{+}$) and $\varepsilon (L^{+})$
are given by

\begin{eqnarray}
\Delta (L^{+}) &=&L\otimes I+I\otimes L-\lambda _1\frac \alpha 4I\otimes
I+\lambda _4(-1)^{\widetilde{N}}\otimes (-1)^{\widetilde{N}}  \nonumber \\
&&-\lambda _1\left( (-1)^{\widetilde{N}}a^{\dagger }\otimes a-(-1)^{-%
\widetilde{N}}a\otimes a^{\dagger }\right) \\
S(L^{+}) &=&L^{+},\quad \varepsilon (L^{+})=\lambda _1\frac \alpha 4+\lambda
_4.
\end{eqnarray}
and that (\ref{hop}) are satisfied. Specialization to include only $L$ in
the above relations, can be obtained by setting $L^{+}\rightarrow L$ and $%
\lambda _4=0$. Relations (\ref{La}) are also preserved by the Hopf algebra
structure (\ref{hop1}), subject to (\ref{con}), while $L^{+}$ is represented
on the Fock space (\ref{foc1}) as

\begin{equation}
L^{+}|n>=\left( \lambda _1(\frac \alpha 4-\frac \beta 2)+\lambda
_4(-1)^{\frac \beta \alpha }\right) (-1)^n|n>  \label{L1}
\end{equation}
Note that $L^{+}$ (and $L$) introduces a $Z_2$ grading on the Fock space.

As mentioned the constraints (\ref{con}) can be widely exploited leading to
various choices of values for $\lambda _i$ in terms of $\alpha ,\beta $. In
the case $\alpha =2$ and $\beta =1,$ $\lambda _1=-\lambda _2,$ $\lambda _3=0$
and $L^{+}=\lambda _1(a^{\dagger }a-N)+\lambda _4(-1)^{\widetilde{N}}.$ Then
(\ref{foc1}) show that on this Fock space $a^{\dagger }a=N,$ $aa^{\dagger
}=N+I$ and $B_{2,1}^{+}$ reduces to the quotient $B/<C>$ (see beginning of
the section) extended with the element $(-1)^{a^{\dagger }a+1/2}$. Also we
can investigate the case where we impose on $B_{\alpha ,\beta }$ (or $%
B_{\alpha ,\beta }^{+}$) the additional relation

\begin{equation}
L^2=\eta I,\quad (\mbox{\rm or (}L^{+})^2=\eta I)  \label{eta}
\end{equation}
with $\eta \in {\bf R},\eta \neq 0$ . From the form of $L$ (\ref{L'}) it can
be easily observed that $L^2$ commutes with all the generators of $B_{\alpha
,\beta }$ (and $B_{\alpha ,\beta }^{+}$) and thus on any faithful
representation it reduces to a multiple of the identity. Also it can be
shown that (\ref{eta}) does not respect the Hopf algebra structure. On a
representation independent way thought using (\ref{L'}) we obtain the
characteristic identity for $C=a^{\dagger }a-\frac \alpha 2N$

\begin{equation}
C\left( C+(\frac \alpha 2-\beta )I\right) +\left( \frac \alpha 4-\frac \beta
2\right) ^2I=\frac \eta {\lambda _1^2}I,  \label{eta1}
\end{equation}
which when solved will give a $L$ as a multiple of the identity. Obviously (%
\ref{eta}) with the choice of $L$ given by (\ref{L'}), {\it is not}
compatible with relations (\ref{La}). The same incompatibility is also true
for the case of $B_{\alpha ,\beta }^{+}$ and $L^{+}$given in (\ref{L''}).
However if we consider the element $\lambda _4(-1)^{\widetilde{N}}=L^{+}$
alone (thus letting $\lambda _1=0$) then (\ref{eta}) can hold (it is just
imposing the requirement $g^2=I$) and the Hopf algebra is preserved.

Finally it is important in what will follow in section 5 to observe that if
we substitute in (\ref{bos}) the generators $N$ obtained from (\ref{L'}) (or
(\ref{L''}) ) as $N=\frac 2\alpha (-\frac 1{\lambda _1}L+a^{\dagger }a+\frac
\alpha 4-\frac \beta 2),$ then (\ref{bos}) becomes

\begin{equation}
\lbrack a,a^{\dagger }]=-\frac 2{\lambda _1}L+\frac \alpha 2I  \label{BH}
\end{equation}
This is true if and only if the values of $\alpha $ and $\beta $ are such
that $L$ and $L^{+}$ contain the monomial $N$, that is when the following
values of the pair ($\alpha ,\beta $) {\it are not} considered: $\alpha
=0,\beta =0$ and $\alpha =0,\beta \neq 0$. Relation (\ref{BH}) shows the
potentiality of $B_{\alpha ,\beta }$(and $B_{\alpha ,\beta }^{+}$) to
accommodate and interchange both commutation and anticommutation relations.
It is interesting to investigate if this relation together with (\ref{La})
serve as an alternative definition of $B_{\alpha ,\beta }$ (or $B_{\alpha
,\beta }^{+}$). This will become clearer in section 5 where (\ref{BH}) and (%
\ref{La}) will be compared with (\ref{vas}).

We can obtain from $B_{\alpha ,\beta }$ a realization of $B(0/1)\simeq
osp(1/2)$ by introducing a {\bf Z}$_2$ grading such that $a$ and $a^{\dagger
}$ are odd and $N$ is even and defining

\begin{eqnarray}
e &=&\mu a^{\dagger },\quad f=\lambda a,\quad h=2N+\frac{2\beta }\alpha I
\label{osp} \\
&&  \nonumber
\end{eqnarray}
provided $\alpha \neq 0$ and $\mu \lambda =\frac 2\alpha $, so that

\begin{equation}
\{e,f\}=h,\quad [h,e]=2e,\quad [h,f]=-2f.
\end{equation}
Then the Hopf algebra structure of $B_{\alpha ,\beta }$ induces a
non-trivial one for $osp(1/2)$ extended by the element $g=(-1)^{\widetilde{N}%
}$ , exactly as in the case of \cite{maj} (but using $B_{\alpha ,\beta }$
instead of the ordinary oscillator algebra) and the $R$-matrix is given by (%
\ref{r0}) below. The Casimir invariant $I_2=-\frac 14e^2f^2-\frac 14ef+\frac
1{16}h^2-\frac 18h$ on the Fock space (\ref{foc1} ) takes the eigenvalue

\begin{equation}
i_2=\frac{\beta ^2}{4\alpha ^2}-\frac \beta {4\alpha }.
\end{equation}
which shows that the representation is irreducible. Similarly we can obtain
a realization of $A_1$ (as a subalgebra of $B(0/1)$ for example) by defining

\begin{equation}
e^{\prime }=\mu ^{\prime }e^2,\quad \mbox{\rm  }f^{\prime }=\lambda ^{\prime
}f^2,\quad h^{\prime }=\frac 12h  \label{A1}
\end{equation}
provided $\mu ^{\prime }\lambda ^{\prime }=-\frac 14$, so that

\[
\lbrack e^{\prime },f^{\prime }]=h^{\prime },\quad [h^{\prime },e^{\prime
}]=2e^{\prime },\quad [h^{\prime },f^{\prime }]=-2f^{\prime }. 
\]
We can also obtain a realization of $sl(2,R)\simeq su(1,1)$ if we set $%
J_0=\frac 12h^{\prime }$, $J_{+}=\frac i{\sqrt{2}}e^{\prime }$ and $%
J_{-}=\frac i{\sqrt{2}}f^{\prime }$. Then the eigenvalues of the $sl(2,R)$
Casimir invariant $C_2=2J_{-}J_{+}-J_0^2-J_0$ on the Fock space (\ref{foc1})
are given by

\begin{equation}
c_n=\frac 12-\frac \beta {2\alpha }-\frac{\beta ^2}{4\alpha ^2}+\frac{2\beta
-\alpha }{8\alpha }(3+(-1)^n).
\end{equation}
which shows that the representation is completely reducible with the two
invariant subspaces corresponding to $n$ being even and $n$ being odd. The
Casimir eigenvalues $c_{even}$ and $c_{odd}$ are 
\[
c_{even}=-\frac 14\frac{\beta ^2}{\alpha ^2}+\frac \beta {2\alpha },\quad
c_{odd}=\frac 14-\frac{\beta ^2}{4\alpha ^2}. 
\]

\subsection{The algebra $\overline{B}_{\sigma ,\tau }$.}

We shall consider now the family of algebras $\overline{B}_{\sigma ,\tau }$
generated by $a,$ $a^{\dagger }$ and $N$ subject to the following relations: 
\begin{eqnarray}
\lbrack a,a^{\dagger }] &=&\sigma N+\tau I,  \nonumber  \label{hos} \\
\lbrack N,a] &=&-a,   \\
\lbrack N,a^{\dagger }] &=&a^{\dagger }  \nonumber
\end{eqnarray}
where $\sigma ,\tau \in {\bf R}$ and for $\sigma =0$ we obtain (\ref{normal}%
). This algebra is a Hopf algebra whose coproduct counit and antipode
satisfying (\ref{hop}) are given by:

\[
\Delta (N)=N\otimes I+I\otimes N+\frac \tau \sigma I\otimes I 
\]

\[
\Delta (a)=a\otimes I+I\otimes a 
\]

\[
\Delta (a^{\dagger })=a^{\dagger }\otimes I+I\otimes a^{\dagger } 
\]

\[
\varepsilon (N)=-\frac \tau \sigma ,\varepsilon (a)=\varepsilon (a^{\dagger
})=0,\varepsilon (I)=1 
\]

\begin{equation}
S(N)=-N-\frac{2\tau }\sigma ,\quad S(a)=-a,\quad S(a^{\dagger })=-a^{\dagger
}  \label{hos1}
\end{equation}
provided $\sigma \neq 0$(thus not allowing a Hopf algebra structure for (\ref
{normal}). An opposite Hopf algebra structure can be easily read off from (%
\ref{hos1}) which also show that an $R$-matrix will turn out to be trivial.
A Fock-type representation, with $a|0>=0,$ $|n>^{\prime }=(a^{\dagger
})^n|0>,$ $N|n>^{\prime }=n|n>^{\prime },n\in Z_{+}$ and $<0|0>=1$, exists
such that

\begin{equation}
a|n>^{\prime }=\overline{[n]}|n-1>^{\prime },\quad \mbox{\rm where}\quad 
\overline{[n]}=\frac{\sigma n(n-1)}2+n\tau .  \label{fo}
\end{equation}
If $\sigma \geq 0,$ $\tau \geq 0$ a unitary representation of $\overline{B}%
_{\sigma ,\tau }$ is provided by:

\[
|n>=\frac 1{(\overline{[n]}!)^{\frac 12}}(a^{\dagger })^n|0> 
\]

\[
a|n>=\overline{[n]}^{\frac 12}|n-1>,\quad a^{\dagger }|n>=\overline{[n+1]}%
^{\frac 12}|n+1>, 
\]

\begin{equation}
\mbox{\rm where}\quad \overline{[n]}!=\prod_{l=1}^n\overline{[l]}\quad \mbox{\rm and}
\quad <n|n^{\prime }>=\delta _{nn^{\prime }}.  \label{fo1}
\end{equation}
Finally an $A_1$ realization can be obtain by setting

\begin{equation}
e=\xi a^{\dagger },\quad \mbox{\rm  }f=\zeta a,\quad h=2N+\frac{2\tau }\sigma I
\mbox{\rm  }
\end{equation}
where $\xi \zeta =-2/\sigma $ and the defining relations of $A_1$ are as
shown below (\ref{A1}).

\section{Deformed Boson algebras $B_{\alpha ,\beta }^q$ and $\overline{B}%
_{\sigma ,\tau }^q$.}

\subsection{The algebra $B_{\alpha ,\beta }^q$.}

We turn now to a $q$-deformation ($q$ generic) of the algebras $B_{\alpha
,\beta }$ and $\overline{B}_{\sigma ,\tau }$. Define $B_{\alpha ,\beta }^q$
as the Lie algebra generated by $a_q$ , $a_q^{\dagger }$ and $N$ subject to
the following relations:

\begin{eqnarray}
\{a_q,a_q^{\dagger }\} &=&\left[ \alpha N+\beta \right] _q,  \nonumber
\label{qbos} \\
\lbrack N,a_q] &=&-a_q,  \nonumber \\
\lbrack N,a_q^{\dagger }] &=&a_q^{\dagger }  \label{qbos1}
\end{eqnarray}
where $\alpha ,$ $\beta \in {\bf R}$ and $[x]_q=(q^x-q^{-x})/(q-q^{-1})$.
This algebra is a Hopf algebra whose coproduct, counit and antipode satisfy (%
\ref{hop}) and are given by:

\[
\Delta (N)=N\otimes I+I\otimes N+\frac \beta \alpha I\otimes I 
\]

\[
\Delta (a)=a_q\otimes q^{\frac{\alpha N+\beta }2}+(-1)^{^{N+\frac \beta
\alpha }}q^{-\frac{\alpha N+\beta }2}\otimes a_q 
\]

\[
\Delta (a_q^{\dagger })=a_q^{\dagger }\otimes q^{\frac{\alpha N+\beta }%
2}+(-1)^{-N-\frac \beta \alpha }q^{-\frac{\alpha N+\beta }2}\otimes
a_q^{\dagger } 
\]

\[
\varepsilon (N)=-\frac \beta \alpha ,\quad\varepsilon (a_q)=\varepsilon
(a_q^{\dagger })=0,\quad\varepsilon (I)=1 
\]

\begin{equation}
S(N)=-N-\frac{2\beta }\alpha ,\quad S(a_q)=-(-1)^{-N-\frac \beta \alpha
}q^{-\frac \alpha 2}a,\quad S(a_q^{\dagger })=a_q^{\dagger }(-1)^{^{N+\frac
\beta \alpha }}q^{\frac \alpha 2}  \label{qhop}
\end{equation}
provided that $\alpha \neq 0$. An opposite Hopf algebra structure also
exists with coproduct $\Delta ^T$ and antipode the inverse $S^{-1}$of $S$
given by :
\begin{equation}
S^{-1}(N)=S(N),\quad S^{-1}(a_q)=(-1)^{-N-\frac \beta \alpha }q^{-\frac
\alpha 2}a_q,\quad S^{-1}(a_q^{\dagger })=-a_q^{\dagger }(-1)^{N+\frac \beta
\alpha }q^{\frac \alpha 2}\mbox{\rm .}
\end{equation}
Similarly to the undeformed case, in order to obtain (\ref{qhop}), we have
to enlarge $B_{\alpha ,\beta }^q$ by adding an invertible element $(-1)^{%
\widetilde{N}}$ which will be treated as an supplementary generator
satisfying relations (\ref{n}) (with $a_q$ and $a_q^{\dagger }$ in the place
of $a$ and $a^{\dagger }$ respectively) and (\ref{n'}). We shall denote this
extended algebra (its universal enveloping algebra) as $B_{\alpha ,\beta
}^{q+}(U(B_{\alpha ,\beta }^{q+}))$.

A Fock-type representation, with $a_q|0>_q=0,$ $|n>_q^{\prime
}=(a_q^{\dagger })^n|0>_q,N|n>_q=n|n>_q,$ $n\in Z_{+}$ and $_q<0|0>_q=1$,
exists such that

\[
a_q|n>_q^{\prime }=(n)_q|n-1>_q^{\prime },\quad (-1)^{\pm N}|n>_q^{\prime
}=(-1)^{\pm n}|n>_q^{\prime }, 
\]

\begin{equation}
\mbox{\rm where }(n)_q=\left( q^{\frac \alpha 2}+q^{-\frac \alpha 2}\right)
^{-1}\left( (-1)^{n+1}[\beta -\frac \alpha 2]_q+[n\alpha +\beta -\frac
\alpha 2]_q\right)  \label{qf}
\end{equation}
Normalizing $|n>_q^{\prime }$ the representation of $B_{\alpha ,\beta }^q$
(and $B_{\alpha ,\beta }^{q+}$) is provided by:

\[
|n>_q=\frac 1{((n)_q!)^{\frac 12}}(a_q^{\dagger })^n|0>_q 
\]

\[
a_q|n>_q=(n)_q^{\frac 12}|n-1>_q,\quad a_q^{\dagger }|n>_q=(n+1)_q^{\frac
12}|n+1>_q, 
\]

\begin{equation}
\mbox{\rm where}\quad (n)_q!=\Pi _{m=1}^n(m)_q\quad \mbox{\rm and}\quad
_q<n|n^{\prime }>_q=\delta _{nn^{\prime }}.  \label{qfoc}
\end{equation}
In the limit $q\rightarrow 1$ we get the Fock space of the undeformed
algebra $B_{\alpha ,\beta }$(or $B_{\alpha ,\beta }^{+}$).

A realization of $osp_{q^{\prime }}(1/2),$ with $q^{\prime }=q^\alpha ,$ can
be obtained by defining

\begin{eqnarray}
e &=&\mu a_q^{\dagger },\quad \mbox{\rm  }f=\lambda a_q,\quad h=N+\frac \beta
\alpha I  \label{qosp} \\
&&  \nonumber
\end{eqnarray}
so that with $\mu \lambda =[\alpha ]_q^{-1}$ the following $osp_{q^{\prime
}}(1/2)$ defining relations are satisfied:

\begin{equation}
\{e,f\}=[h]_{q^{\prime }},\quad [h,e]=e,\quad [h,f]=-f.
\end{equation}
Finally an {\bf $R$-}matrix exists and is given by

\begin{equation}
R=R_0q^{\alpha \widetilde{N}\otimes \widetilde{N}}\sum_{l=0}^\infty
(q-q^{-1})^lq^{-\frac \alpha 4l(l+1)}\frac{(-1)^{\frac 14l(l-1)}}{[l]_x!}%
q^{\frac \alpha 2l\widetilde{N}^{}}(-1)^{l\widetilde{N}}{\LARGE 
(a_q^{\dagger })^{{\it l}}\otimes }q^{-\frac \alpha 2l\widetilde{N}^{}}
{\LARGE a_q^{{\it l}}}  \label{r11}
\end{equation}
where $\widetilde{N}=N+{\beta }/{\alpha }$ , $x=(-q^{-\alpha })^{1/2}$ and 
\begin{equation}
R_0=\frac 12(I\otimes I+I\otimes (-1)^{\widetilde{N}}+(-1)^{\widetilde{N}%
}\otimes I-(-1)^{\widetilde{N}}\otimes (-1)^{\widetilde{N}}).  \label{r0}
\end{equation}
provided that we demand that $(-1)^{2\widetilde{N}}=I$. (\ref{r1}) has been
calculated using quantum double technics similar to \cite{dt}. It is
important to mention that (\ref{r1}) is exactly the same as the one appeared
in relation (49) in \cite{maj} with $q\rightarrow q^\alpha ,$ provided we do
the following identifications with the generators of the bosonization of $%
osp_{q^\alpha }(1/2)$ : $J_z=\frac 12\widetilde{N},$ $V_{+}=ka_q^{\dagger },$
$V_{-}=ta_q,$ $g=(-1)^{\widetilde{N}},$ and $kt=-[4]_{q^\alpha }^{-1}[\alpha
]_q^{-1}$. Then we can argue that $B_{\alpha ,\beta }^{q+}$ is the spectrum
generating quantum group for the ordinary $q$-deformed harmonic oscillator
defined by the relations $a_qa_q^{\dagger }-q^{\pm }a_q^{\dagger }a_q=q^{\mp
N}$ and the last two of (44), which as will be shown in \cite{pa} is a
subalgebra of $B_{\alpha ,\beta }^{q+}.$ At the limit $q\rightarrow 1,$ $%
R\rightarrow R_0$ which is the $R-$ matrix of the undeformed $B_{\alpha
,\beta }^{+}$.

\subsection{The algebra $\overline{B}_{\sigma ,\tau }^q$.}

Similarly define $\overline{B}_{\sigma ,\tau }^q$ as the Lie algebra
generated by $a_q$ , $a_q^{\dagger }$ and $N$ subject to the following
relations:

\begin{eqnarray}
\lbrack a_q,a_q^{\dagger }] &=&\left[ \sigma N+\tau I\right] _q,  \nonumber
\label{qbos11} \\
\lbrack N,a_q] &=&-a_q,  \label{qhos} \\
\lbrack N,a_q^{\dagger }] &=&a_q^{\dagger }  \nonumber
\end{eqnarray}
where $\sigma ,$ $\tau \in {\bf R}$.

This algebra is a Hopf algebra whose coproduct, counit and antipode satisfy (%
\ref{hop}) and are given by:

\[
\Delta (N)=N\otimes I+I\otimes N+\frac \tau \sigma I\otimes I,\quad \lambda
\neq 0 
\]

\[
\Delta (a_q)=a_q\otimes q^{\frac{\sigma N+\tau }2}+q^{-\frac{\sigma N+\tau }%
2}\otimes a_q 
\]

\[
\Delta (a_q^{\dagger })=a_q^{\dagger }\otimes q^{\frac{\sigma N+\tau I}%
2}+q^{-\frac{\sigma N+\tau }2}\otimes a_q^{\dagger } 
\]

\[
\varepsilon (N)=-\frac \tau \sigma ,\quad \varepsilon (a_q)=\varepsilon
(a_q^{\dagger })=0,\quad \varepsilon (I)=1 
\]

\begin{equation}
S(N)=-N-\frac{2\tau }\sigma ,\quad S(a_q)=-q^{-\frac \sigma 2}a_q,\quad
S(a_q^{\dagger })=-q^{\frac \sigma 2}a_q^{\dagger }.  \label{qhop1}
\end{equation}
An opposite Hopf algebra structure exists whose antipode $S^{-1}$ is given by

\[
S^{-1}(N)=S(N),\quad S^{-1}(a_q)=-q^{\frac \sigma 2}a_q,\quad
S^{-1}(a_q^{\dagger })=-q^{-\frac \sigma 2}a_q^{\dagger } 
\]
A Fock-type representation, with $a_q|0>_q=0,$ $|n>_q^{\prime
}=(a_q^{\dagger })^n|0>_q$ $N|n>_q^{\prime }=n|n>_q^{\prime },$ $n\in Z_{+}$%
, and $_q<0|0>_q=1$, exists such that normalizing $|n>_q^{\prime }$we get

\[
|n>_q=\frac 1{(\overline{(n)}_q!)^{\frac 12}}(a_q^{\dagger })^n|0>_q 
\]

\begin{eqnarray*}
a_q|n &>&_q=\overline{(n)}_q^{\frac 12}|n-1>_q,\quad a_q^{\dagger }|n>_q=%
\overline{(n+1)}_q^{\frac 12}|n+1>_q \\
&&
\end{eqnarray*}

\[
\overline{(n)}_q=\left[ \frac \sigma 2\right] _q^{-1}\left( [\frac{\sigma n}%
2]_q[\frac{\sigma (n-1)}2+\tau ]_q\right) 
\]

\begin{equation}
\overline{(n)}_q!=\Pi _{m=1}^n\overline{(m)}_q\quad \mbox{\rm  and}\quad \mbox{\rm  }%
_q<n|n^{\prime }>_q=\delta _{nn^{\prime }}  \label{qfo}
\end{equation}
In the limit $q\rightarrow 1$ we get the Fock space of the undeformed
algebra $\overline{B}_{\sigma ,\tau }$. A $sl_{q^{\sigma /2}}(2)$
realization is provided by the following identifications

\begin{equation}
e=\xi a_q^{\dagger },\quad \mbox{\rm  }f=\zeta a_q,\quad \mbox{\rm  }h=2N+\frac{%
2\tau }\sigma I\quad
\end{equation}
where $\xi \zeta =-[\frac \sigma 2]_q$, so that the defining $sl_{q^{\sigma
/2}}(2)$ relations below are satisfied:

\begin{equation}
\lbrack e,f]=\left[ h\right] _{q^{\sigma /2}},\quad [h,e]=2e,\quad [h,f]=-2f.
\end{equation}
An {\bf $R$-}matrix also exists, and is given by

\begin{equation}
R=q^{\sigma \widetilde{N}\otimes \widetilde{N}}\sum_{l=0}^\infty q^{\frac
\sigma 4l(l+1)^{}}\frac{(-1)^l}{[l]_x!}(q-q^{-1})^lq^{\frac \sigma 4l
\widetilde{N}^{}}{\LARGE (a_q^{\dagger })^{{\it l}}\otimes }q^{-\frac \sigma
4l\widetilde{N}^{}}{\LARGE a_q^{{\it l}}}  \label{r22}
\end{equation}
where $\widetilde{N}=N+\frac \tau \sigma $ and $x=(q^{{\sigma }/2})$. As in
the case of (\ref{r1}), (\ref{r2}) has also been obtained using quantum
double technic.

\section{The generalized `$\nu -$deformed' Heisenberg algebra H$_{\delta
,\nu }$.}

We shall generalize now the so-called `deformed' Heisenberg algebra of
Vassiliev \cite{vas1}. This is defined as the algebra $H_{\delta ,\nu }$
generated by $b$, $b^{\dagger }$ and $K$ subject to the following relations:

\begin{eqnarray}
\lbrack b,b^{\dagger }] &=&\delta I+\nu K  \nonumber  \label{vas} \\
\{K,b\} &=&\{K,b^{\dagger }\}=0  \label{vasi}
\end{eqnarray}
where $\delta ,\nu \in {\bf R}$. If we impose the additional requirements
that $K^2=I$ then with $\delta =1$ we obtain that of \cite{vas1}, used for
example in \cite{mac2}, \cite{pl1}, \cite{pl2}.

A Fock-type representation ( a generalization of that appearing in \cite{pl1}%
, \cite{pl2} ), with $b|0>=0$, $|n>^{\prime }=(b^{\dagger })^m|0>,$ $m\in 
{\bf Z}_{+}$ and $<0|0>=1$, exists such that

\begin{equation}
b|n>^{\prime }=[m]|m-1>^{\prime },\quad \mbox{\rm where}\quad [m]=\delta m+\frac{%
\nu -\delta +1}2(1+(-1)^{m+1})  \label{Kf'}
\end{equation}
If $\nu >0$ and $\delta >0$ a unitary representation of $B_{\alpha ,\beta }$
is provided by:

\[
|m>=\frac 1{([m]!)^{\frac 12}}(b^{\dagger })^m|0>, 
\]

\[
b|m>=[m]^{\frac 12}|m-1>,\quad b^{\dagger }|m>=[m+1]^{\frac 12}|m+1>, 
\]

\[
K|m>=\frac{\nu -\delta +1}\nu (-1)^m|m>, 
\]

\begin{equation}
\lbrack m]!=\Pi _{l=1}^m[l],\quad <m|m^{\prime }>=\delta _{mm^{\prime }}.
\label{Kf}
\end{equation}
The striking similarity of the Fock spaces (\ref{foc1}) and (\ref{Kf}) is
not accidental. As we shall just show under certain conditions we can obtain 
$B_{\alpha ,\beta }$ from $H_{\delta ,\nu }$ and vice versa, not only on the
above Fock spaces but as abstract algebras. $H_{\delta ,\nu }$ can be
extended so that the resulting algebra will possess a Hopf algebra
structure. There exists in the enveloping algebra $U(H_{\delta ,\nu })$ of $%
H_{\delta ,\nu }$ , an element $M$ given by

\begin{equation}
M=\mu _1b^{\dagger }b+\mu _2K+\rho I  \label{Ma}
\end{equation}
and satisfying

\begin{equation}
\lbrack M,b]=-b,\quad [M,b^{\dagger }]=b^{\dagger },  \label{Mc}
\end{equation}
where $\mu _i,\rho \in {\bf R}$ $(i=1,2)$ provided that the following
constraint is satisfied

\begin{equation}
\mu _1\delta I+(2\mu _2-\nu \mu _1)K=I.  \label{con2}
\end{equation}
This suggests that, since $K$ {\it should not} be a multiple of the
identity, (as this contradicts the second of (\ref{vas})) necessarily $%
\delta \neq 0$ which leads to $\mu _1=1/\delta $ and $2\mu _2-\nu \mu _1=0$.
Consequently with the above constraints (\ref{Ma}) now becomes

\begin{equation}
M=\frac 1\delta b^{\dagger }b+\frac \nu {2\delta }K+\rho I,  \label{Ma'}
\end{equation}

\begin{equation}
\mbox{\rm and}\quad M|m>=(m+\frac{\nu -\delta +1}{2\delta }+\rho )|m>
\end{equation}
which by choosing $\rho =-\frac{\nu -\delta +1}{2\delta }$, $M|m>=m|m>$. The
choice of $M$ given by (\ref{Ma'}) is obviously not the most general
possible but it is the unique non-zero combination of lowest order monomials
of generators of $H_{\delta ,\nu }$ that satisfy (\ref{Mc}) while such an
element does not exist if $\delta =0$ (these considerations can be inferred
by writing $M^{\prime }=C_{lmn}K^l(b^{\dagger })^mb^n$ , $l$, $m$, $n$ $\in 
{\bf Z}_{+}$ , $C_{lmn}\in {\bf R}$ , and demanding that (\ref{Mc}) are
satisfied together with $[M^{\prime },K]=0)$. Now we are in a position to
demonstrate the similarities between $B_{\alpha ,\beta }$ and $H_{\delta
,\nu }$. Solving (\ref{Ma'}) with respect to $K$ and substituting in the
first of (\ref{vas}) we obtain

\begin{equation}
\{b,b^{\dagger }\}=2\delta M+\delta (1-2\rho )I  \label{com}
\end{equation}
so that together with (\ref{Mc}) $H_{\delta ,\nu }$ takes the form of the
defining relations of $B_{\alpha ,\beta },$ (\ref{bos}) by setting

\begin{equation}
\alpha =2\delta ,\quad \beta =\delta -2\rho \delta  \label{hb}
\end{equation}
and where $M$ is replaced by $N$. Note that for the case where we choose $%
\rho =-\frac{\nu -\delta +1}{2\delta },$ $H_{\delta ,\nu }$ takes the form
of $B_{2\delta ,\nu +1}$. This process can also be carried out the opposite
direction, as (\ref{La}) and (\ref{BH}) suggests, by setting in (\ref{BH}) 
\begin{equation}
\delta =\alpha /2,\quad \nu =-2/\lambda _1  \label{bh}
\end{equation}
and where $L$ is replaced by $K$. It is easy to observe that, as $\lambda
_1\neq 0$, (\ref{bh}) shows that $B_{\alpha ,\beta }$ can not be mapped to a 
$H_{\delta ,0}$ -form. Also $B_{0,\beta }$ can not be mapped to a $H_{\delta
,\nu }$ -form at all, since the appropriate $L$ fails to exist (no monomial $%
N$ is present in $L$ even if we perform a thorough search for a more general 
$L$ in $U(B_{\alpha ,\beta })$). Also (\ref{hb}) shows that a $B_{\alpha
,\beta }$ -form of $H_{0,\nu }$ fails to exist since $M$ can not be defined
and for $H_{\delta ,0}$ the appropriate $M$ does not exist (no monomial $K$
is present in $M$) thus also not allowing a $B_{\alpha ,\beta }$ -form.
Consequently provided that we keep away from the values $\alpha =$ $\delta =$
$\nu =0$ we can always obtain a $H_{\delta ,\nu }$ -form of $B_{\alpha
,\beta }$ and vice-versa. Relations (\ref{hb}) and (\ref{bh}) also imply
that $\rho $ and $\beta $ can have arbitrary values. However observation of
the Fock spaces (\ref{foc1}) and (\ref{Kf}) and comparison of the action of $%
K$, $M$, $L$ and $N$ on them, shows that with the identifications (\ref{hb}%
), (\ref{bh}), $\rho =-\frac{\nu -\delta +1}{2\delta }$ , $\beta =\nu +1$
not only these spaces are equivalent but also $H_{\delta ,\nu }$ and $%
B_{\alpha ,\beta }$ are isomorphic with $K\equiv L$, $N\equiv M$, $b\equiv a$%
, $b^{\dagger }\equiv a^{\dagger }$.

It can be checked that the following maps $\varphi :$ $B_{\alpha ,\beta
}\rightarrow H_{\delta ,\nu }$ and $\varphi ^{\prime }:$ $H_{\delta ,\nu
}\rightarrow B_{\alpha ,\beta }$ defined by:

\begin{equation}
\varphi (a)=b,\quad \varphi (a^{\dagger })=b^{\dagger },\quad \varphi
(N)=\frac 2\alpha b^{\dagger }b+\frac \nu \alpha K+\frac{\delta -\beta }%
\alpha ,\quad \alpha \neq 0,  \label{is1}
\end{equation}

\begin{equation}
\varphi ^{\prime }(b)=a,\quad \varphi ^{\prime }(b^{\dagger })=a^{\dagger
},\quad \varphi ^{\prime }(K)=-\frac 2\nu a^{\dagger }a+\frac \alpha \nu N+%
\frac{\beta -\delta }\nu ,\quad \nu \neq 0  \label{is2}
\end{equation}
are homomorphisms {\it if and only if} $\alpha =2\delta $ . Moreover $%
\varphi ^{\prime }=\varphi ^{-1}$ and $\varphi $ becomes an isomorphism $%
H_{\delta ,\nu }$ $\simeq B_{\alpha ,\beta }$ {\it provided} that both $%
\alpha \neq 0$ {\it and} $\nu \neq 0$. $\varphi $ and $\varphi ^{\prime }$
can be thought as defining families of maps where each member is
parametrized by $\alpha ,\beta ,\nu $ and $\delta ,\beta ,\nu $ respectively
and we can formally write $\varphi \equiv \varphi _{\alpha ,\beta ,\nu }$
and $\varphi ^{\prime }\equiv \varphi _{\delta ,\nu ,\beta }^{\prime }$. So
for example $B_{2,1}$is mapped via $\varphi _{2,1,\nu }$ to $H_{1,\nu }$ ($%
\nu $ a fixed chosen number) and $H_{1,\nu }$ is mapped via $\varphi
_{2,1,\nu }^{-1}=\varphi _{1,\nu ,1}$ back to $B_{2,1}$. Finally it can be
checked that

\begin{equation}
\varphi (L)=-\frac{\lambda _1\nu }2K,\quad \varphi ^{\prime }(M)=N+\frac{%
\beta -\delta }{2\delta }+\rho .
\end{equation}

With the existence of $M$ of the form (\ref{Ma'}) we can enlarge $H_{\delta
,\nu }$ by adding the invertible element $(-1)^M$, as we did in the case of $%
B_{\alpha ,\beta }$. We shall denote this enlarged algebra by $H_{\delta
,\nu }^{+}$ and the relations that $(-1)^M$ has to satisfy are given by:

\begin{eqnarray}
\{(-1)^{\pm M},a\}=\{(-1)^{\pm M},a^{\dagger }\}=0,\quad [(-1)^{\pm
M},K]=0,\quad  \label{m}
\end{eqnarray}
and on (\ref{Kf}) will be represented as $(-1)^M|m>=(-1)^{m+\frac{\nu
-\delta +1}{2\delta }+\rho }|m>$. $H_{\delta ,\nu }^{+}$ can obviously be
treated in the spirit of \cite{maj} as was done with $B_{\alpha ,\beta }^{+}$
with the element $g$ of \cite{maj} being $g=(-1)^{\widetilde{M}},$ where $%
\widetilde{M}=M-\rho +\frac 12$. Thus $H_{\delta ,\nu }^{+}$ can be
considered as a spectrum generating algebra of the ordinary oscillator
algebra. Moreover the isomorphism can also be extended such that $H_{\delta
,\nu }^{+}$ $\simeq B_{\alpha ,\beta }^{+}$ by defining $\varphi ((-1)^{%
\widetilde{N}})=(-1)^{\widetilde{M}},\varphi ^{\prime }((-1)^{\widetilde{M}%
})=(-1)^{\widetilde{N}}$.

Then a Hopf algebra structure for $H_{\delta ,\nu }^{+}$ is given by:

\begin{eqnarray*}
\Delta (K) &=&K\otimes I+I\otimes K+\frac \delta \nu I\otimes I \\
&&-\frac 2\nu (-1)^{-M+\rho -\frac 12}b\otimes b^{\dagger }+\frac 2\nu
(-1)^{M-\rho +\frac 12}b^{\dagger }\otimes b
\end{eqnarray*}

\[
\Delta (b)=b\otimes I+(-1)^{M-\rho +\frac 12}\otimes b 
\]

\[
\Delta (b^{\dagger })=b^{\dagger }\otimes I+(-1)^{-M+\rho -\frac 12}\otimes
b^{\dagger } 
\]

\[
\varepsilon (K)=-\frac \delta \nu ,\quad \varepsilon (b)=\varepsilon
(b^{\dagger })=0,\quad \varepsilon (I)=1 
\]

\begin{equation}
S(K)=K,\quad S(b)=-(-1)^{-M+\rho -\frac 12}b,\quad S(b^{\dagger
})=b^{\dagger }(-1)^{M-\rho +\frac 12}  \label{Hhop}
\end{equation}

\begin{eqnarray}
\Delta ((-1)^{\pm M}) &=&(-1)^{\pm (\frac 12-\rho )}(-1)^{\pm M}\otimes
(-1)^{\pm M},  \nonumber \\
S((-1)^{\pm M}) &=&(-1)^{\mp M\pm (2\rho -1)},\quad \varepsilon ((-1)^{\pm
M})=(-1)^{\pm (\rho -\frac 12)}  \label{m'}
\end{eqnarray}
provided that $\nu ,\delta \neq 0$. As in the case of $B_{\alpha ,\beta
}^{+} $ it is not necessary to impose at this stage the condition $(-1)^{2%
\widetilde{M}}=I$ (which implies that $(-1)^{2M}=(-1)^{2\rho -1}$ ). The
form of $\Delta (M)$, $\varepsilon (M)$ and $S(M)$ is given by

\[
\Delta (M)=M\otimes I+I\otimes M+(\frac 12-\rho )I\otimes I, 
\]

\begin{equation}
\varepsilon (M)=\rho -\frac 12,\quad S(M)=-M+2\rho -1.  \label{Mhop}
\end{equation}
It can be checked that using $\varphi $ and $\varphi ^{\prime }$ we can show
that the above mentioned isomorphism carries to the Hopf algebra structures
of $B_{\alpha ,\beta }^{+}$ and $H_{\delta ,\nu }^{+}$ too. An opposite Hopf
algebra structure also exists with an antipode the inverse of the one given
above:

\begin{equation}
.S^{-1}(K)=K,\quad S^{-1}(b)=(-1)^{-\widetilde{M}}b,\quad S^{-1}(b^{\dagger
})=-b^{\dagger }(-1)^{\widetilde{M}},\quad S^{-1}((-1)^{\widetilde{M}%
})=(-1)^{\widetilde{M}}
\end{equation}

Finally we can obtain a realization of $osp(1/2)$ by defining

\begin{equation}
e=\mu b^{\dagger },\quad f=\lambda b,\quad h=\frac \nu \delta K+\frac
2\delta b^{\dagger }b+I
\end{equation}
provided that $\mu \lambda =\frac 1\delta $, while for $A_1$ (as a
subalgebra of $osp(1/2)$) by defining

\begin{equation}
e^{\prime }=\mu ^{\prime }(b^{\dagger })^2,\quad f=\lambda ^{\prime
}b^2,\quad h^{\prime }=\frac \nu {2\delta }K+\frac 1\delta b^{\dagger
}b+\frac 12I
\end{equation}
provided $\mu ^{\prime }\lambda ^{\prime }=-\frac 1{4\delta ^2}$.
Implementing the Hopf structure of $H_{\delta ,\nu }$ we can obtain a Hopf
structure for the bosonization of $osp(1/2)$ as was the case for $B_{\alpha
,\beta }^{+}$ . In particular it is expected that an $R$ matrix for $%
H_{\delta ,\nu }$ will be of the form of (\ref{r0}) with $\widetilde{M}$ in
the place of $\widetilde{N}$ and provided we also demand that $(-1)^{2%
\widetilde{M}}=I$.

\section{Conclusion}

In this paper we considered the generalized boson algebras $B_{\alpha ,\beta
}$ and their $q$-deformed versions $B_{\alpha ,\beta }^q$ which when
enlarged by the element $(-1)^N$ it was shown to admit a quasitriangular
Hopf algebra structure. In particular this structure revealed the property
that $B_{\alpha ,\beta }^{+}$ and $B_{\alpha ,\beta }^{q+}$ can be treated
as spectrum generating quantum groups for the undeformed and $q$-deformed
bosons respectively. An important point though is the isomorphism of $%
B_{\alpha ,\beta }$ and $H_{\delta ,\nu }$ (and their respective
enlargements $B_{\alpha ,\beta }^{+}$ and $H_{\delta ,\nu }^{+}$)
demonstrated in section 5 which carries over to their Hopf algebra
structure. To our knowledge it is the first time that the Calogero-Vasiliev $%
\nu $-deformed Heisenberg algebra $H_{1,\nu }$ , slightly modified (i.e. $%
K^2\neq I$ ), can be formulated as a Hopf algebra. Moreover it is expected
that there should exist a $q$-deformation, $H_{\delta ,\nu }^q$ other than
the one of \cite{macfar} or \cite{mac2} which may admit a Hopf structure,
giving a two-parameter deformation of the Heisenberg algebra and possibly
not being isomorphic with $B_{\alpha ,\beta }^q$ , thus giving rise to a new 
$R$-matrix. Consequences of these Hopf-type boson algebras on physical
models such as the Calogero-Sutherland models, supersymmetric quantum
mechanics, anyonic systems (whose references mentioned in the introduction)
or on radial problems, BRST symmetry \cite{jar1}, are under investigation.
It is anticipated that the Hopf algebra structure, and especially the
quasitriangular nature of these algebras, might reveal interesting
connections with the integrability of the above physical systems.

Another important aspect of these models is their relations with existing
ones. In work under completion, \cite{pa}, we investigate the various
subalgebras of these undeformed and deformed models using the powerful tool
of the fixed point set of the adjoint action of a Hopf algebra. It is shown
how already known undeformed and $q$-deformed boson algebras appear as fixed
point subalgebras or as appropriate quotients. It is at this point that the
role of the Cuntz algebra is also investigated. Finally the more natural
generalization of the ordinary boson algebra, $\overline{B}_{\sigma ,\tau }$%
, was considered together with its $q$-deformation, and both were shown to
be quasitriangular Hopf algebras too.

Braid group representations and possible link invariants for all of the
proposed models of deformed and undeformed bosons are worth investigating,
while the $osp(1/2)$ and $A_1$ realizations obtained point towards
realizations of higher rank algebras and superalgebras which will also allow
for the construction of families of infinite dimensional representations
when the above Fock-spaces are generalized.

Finally one should comment on the implications of the generalized boson
algebras in particular $B_{\alpha ,\beta },$ and $B_{\alpha ,\beta }^q,$ for
quantum statistics. As the usual oscillator algebra does not possess a Hopf
algebra structure, it is difficult to characterize the multiparticle
Hamiltonian. However, in our case by generalizing to a many-particle system $%
B_{\alpha ,\beta }^i$ $i=1,2,...($ and in particular taking $\alpha =2,\beta
=1$ which on the Fock space will give $a^{i\dagger }a^i=N^i,$ $%
a^ia^{i\dagger }=N^i+I$) the total Hamiltonian (taken to be proportional to $%
\alpha \widetilde{N}=\{a,a^{\dagger }\}$) has a very natural interpretation
as being proportional to $\Delta ^{(n)}(\widetilde{N}),$ the $n$-fold
coproduct of the one particle Hamiltonian. Perhaps of most interest are
implications for the quantum statistics of $B_{\alpha ,\beta }^q$. In this
case the existence of the coproduct of $\widetilde{N}$ implies various
logical possibilities for the multiparticle Hamiltonian, which may be more
acceptable than the obvious (but arbitrary) choice $\propto \sum_i[N_i]_q$
which has no justification in terms of a Hopf structure. Non-local effects
will likely emerge from such choices, which may play a crucial role in
modifying the partition functions and statistics of the system.

{\bf Acknowledgments}

The authors would like to thank P. E. T. J\o rgensen, A. J. Bracken, D. S.
McAnally and R. Zhang for their sincere interest, support and fruitful
comments during the completion of this work parts of which were reported in
the 12th Australian Institute of Physics Congress, Hobart July 1-5, 1996 and
in the 3rd International Conference on Functional Analysis and Approximation
Theory, September 23-28, 1996 Acquafredda di Maratea, Potenza, Italy \cite
{ours}.

\end{document}